\documentclass{article}
\usepackage{amsthm}
\usepackage{amsfonts}
\usepackage{amssymb}
\usepackage{amsgen}
\usepackage{amsmath}
\usepackage{amsopn}
\usepackage{verbatim}
\usepackage{xypic}
\usepackage{xspace}
\usepackage{multicol}
\usepackage{makeidx}
\usepackage{eepic}
\usepackage{upref}

\theoremstyle{plain}
\newtheorem{thm}{Theorem}[section]
\newtheorem{lem}[thm]{Lemma}
\newtheorem{prop}[thm]{Proposition}
\newtheorem{cor}[thm]{Corollary}

\theoremstyle{definition}

\newtheorem{defn}[thm]{Definition}
\newtheorem{rem}[thm]{Remark}

\newtheorem{eg}[thm]{Example}

\newcommand{\barp}{{\bar{p}}}

\newcommand{\ZN}{{\mathbb Z_{\ge 0}}}

\DeclareMathOperator{\re}{{Re}}

\newcommand{\ga}{{\alpha}}

\begin{document}

\title{Bernoulli Numbers, Wolstenholme's Theorem, and
$p^5$ Variations of Lucas' Theorem }
\author{Jianqiang Zhao\footnote{Partially supported by
NSF grant DMS0139813}}
\date{}
\maketitle
\begin{center}
{\large Department of Mathematics, Eckerd College, FL 33711}
\end{center}

\textbf{Abstract.} In this note we shall improve some congruences
of G.S.~Kazandzidis and D.F.~Bailey to higher prime power moduli,
by studying the relation between irregular pairs of the form
$(p,p-3)$ and refined version of Wolstenholme's theorem.

\section{Introduction}
Let $H_1(n)$ be the $n$th partial sum of the harmonic series. It
is a classical result commonly attributed to Wolstenholme
\cite[p.89]{HW} that for any prime $p\ge 5$
\begin{equation}\label{equ:wols} H_1(p-1):=
\sum_{k=1}^{p-1} \frac{1}{k}\equiv 0 \pmod{p^2}.
\end{equation}
It's also known \cite{ELeh} that $H_1(p-1)\equiv 0 \pmod{p^3}$ if
and only if $(p,p-3)$ is an irregular pair, namely, $p$ divides
the numerator of $B_{p-3}$. Here we define the Bernoulli numbers
$B_k$ by the Maclaurin series
$$\frac{x}{e^x-1}=\sum_{k=0}^\infty B_k\frac{x^k}{k!}.$$

There is another important equivalent statement of Wolstenholme's
Theorem by using combinatorics. D.F. Bailey \cite{Bailey}
generalizes it to the following form.
\begin{thm}\label{b1}
{\em (\cite[Theorem~4]{Bailey})} Let $n$ and $r$ be non-negative
integers and $p\ge 5$ be a prime. Then
$${np \choose rp} \equiv {n \choose r} \pmod{p^3},$$
where we set ${n \choose r}=0$ if $n<r$.
\end{thm}
He further obtains the following variation of Lucas' Theorem.
\begin{thm}  \label{b2} {\em (\cite[Theorem~5]{Bailey})}
Let $N$, $R$, $n$ and $r$ be non-negative integers and $p\ge 5$ be
a prime. Suppose $n,r<p$. Then
$${Np^3+n \choose Rp^3+r} \equiv {N \choose R} {n \choose r} \pmod{p^3}.$$
\end{thm}

In late 1960's G.~S.~Kazandzidis worked on similar congruences.
Define for any integer $n$ and any positive integer $r$
$$\overline{{n\choose 0}}=1, \quad
 \overline{{n\choose r}}=\frac{n(n+1)\cdots(n+r-1)}{r!}.$$
Among many results he obtained in \cite{K1,K2} the followings are
particular relevant to our study
\begin{thm}  \label{K1} {\em (\cite[$2^{\ast\ast}$ on p.~10]{K1})}
Let $n$ be any integer and $r$ be any positive integers and $p\ge
3$ be a prime. Then
\begin{equation}\label{equ:K1}
\overline{{np \choose rp}}\Big{/} \overline{{n \choose r}}\equiv
 \left\{ \aligned & 1-p^2nr(n+r) &\pmod{ p^3 } & \quad &\text{ if }p=3\ \\
&1 &\pmod{ p^3 }&\ &\text{if }p>3,\endaligned \right.
\end{equation}
and
\begin{equation}\label{equ:K2}
{np \choose rp}\Big{/}  {n \choose r} \equiv
 \left\{ \aligned & 1-p^2nr(n-r)   &\pmod{ p^3 } & \quad &\text{ if }p=3\ \\
&1 &\pmod{ p^3 }&\  &\text{if }p>3.\endaligned \right.
\end{equation}
Here \eqref{equ:K1} and \eqref{equ:K2} are equivalent.
\end{thm}

In this short note we will refine the above results for primes
$p>5$ by using higher prime power modulus (see
Theorem~\ref{thm:1}. This is best possible in the sense that the
result would be wrong if we allowed $p=5$. Note that in \cite{K2}
Kazandzidis obtains an improved version of his congruences of
Theorem \ref{K1} by  replacing the modulus by $p^3 \barp
\{nr(n-r)\}$, where $\barp\{N\}$ denotes the highest power of the
prime $p$ that divides $N$. This improvement does not follow from
our result in this paper. However, it does not imply ours either.

In the last section of this paper we provide an interesting
congruence involving Bernoulli numbers of the form $B_{p-3}$ for
odd primes $p$. The author wishes to thank the referee for
pointing out Kazandzidis' work \cite{K1,K2} and many other
valuable comments which makes this note more readable.

\section{Preliminaries and Some Notation}
Define the Euler-Zagier multiple zeta functions of depth $d$ by
\begin{equation}\label{zeta}
\zeta(s_1,\dots, s_d)=\sum_{0<k_ 1<\dots<k_d} k_1^{-s_1}\cdots
k_d^{-s_d}
\end{equation}
for complex variables $s_1,\dots, s_d$ satisfying
$\re(s_j)+\dots+\re(s_d)>d-j+1$ for all $j=1,\dots,d$. The special
values of multiple zeta functions at positive integers have
significant arithmetic and algebraic meanings, whose defining
series \eqref{zeta} will be called {\em MZV series}, and whose
$n$th partial sum is
\begin{equation}\label{equ:defn}
H(s_1,\dots,s_d;n):=\sum_{1\le k_ 1<\dots<k_d\le n}
k_1^{-s_1}\cdots k_d^{-s_d},\quad n\in \ZN.
\end{equation}
Note that partial sums exist even for divergent MVZ
$\zeta(\dots,1)$ such as the harmonic series $\zeta(1)$. When an
ordered set $(e_1,\dots,e_t)$ is repeated $d$ times we abbreviate
it as $\{e_1,\dots, e_t\}^d$. From the definitions \eqref{zeta}
and \eqref{equ:defn} one derives easily the so called shuffle
relations. For example
$$\zeta(s)\zeta(t)=\zeta(t,s)+\zeta(t+s)+\zeta(s,t)$$
because
$$\sum_{k>0}\cdot \sum_{l>0}=\sum_{k>l>0} +
 \sum_{k=l>0}+\sum_{0<k<l}.$$
Similarly, one has
\begin{equation}\label{shuf}
H(s;n)H(t;n)=H(t,s;n)+H(t+s;n)+H(s,t;n).
\end{equation}

Recall that Stirling numbers $S(n,j)$ of the first kind are
defined by the expansion
\begin{equation}\label{fx}
 f_n(x)=x(x-1)(x-2)\cdots(x-n+1)=\sum_{j=1}^n (-1)^{n-j} S(n,j) x^j.
\end{equation}
These numbers are related to the partial sums of nested harmonic
series:
\begin{equation}\label{Stirling}
S(n,j) = (n-1)!  H(\{1\}^{j-1}; n-1) , \text{ for } j=1,\cdots, n.
\end{equation}
For example, $S(n,n)=1$, $S(n,n-1)=n(n-1)/2$, and $S(n,1)=(n-1)!$.
In particular, if $n=p$ is a prime we then have
$$f_p(x)=(p-1)!x(1-H(1;p-1)x+H(1,1;p-1)x^2-\cdots+x^{p-1}).$$
Comparing $f_p(p)=p!$ we recover the Wolstenhomle's Theorem.

One last thing we need in this note is the following
generalization of \eqref{equ:wols}.
\begin{lem}\label{genWols}  {\em (\cite[Theorem~2.13]{wols})}
Let $s$ and $d$ be two positive integers.  Let $p$ be an odd prime
such that $p\ge sd+3$. Then
$$H(\{s\}^d; p-1) \equiv
\begin{cases}
0 \pmod{p^2} \quad&\text{if } 2\nmid sd,\\
0 \pmod{p} \quad&\text{if } 2\mid sd.
\end{cases}$$
\end{lem}

\section{Main Results}
Our first result improves on Theorem~\ref{b1} of Bailey and
Theorem~\ref{K1} of Kazandzidis simultaneously for all primes
greater than 5.

\begin{defn} For any prime $p\ge 5$ by Wostenholme's theorem
\eqref{equ:wols} we define $w_p<p^2$ to be the unique non-negative
integer such that $w_p\equiv H_1(p-1)/p^2\pmod{p^2}$. It is a
well-known fact that (see for eg. \cite{Gl})
 \begin{equation}\label{Gl}
w_p\equiv -\frac13 B_{p-3}\pmod{p}.
\end{equation}
\end{defn}
\begin{thm}\label{thm:1}
Let $n$ and $r$ be non-negative integers and $p\ge 7$ be a prime.
Then
\begin{equation}\label{main}
{np \choose rp}\Big{/} {n \choose r}\equiv
 1+w_pnr(n-r)p^3 \pmod{p^5 }.
\end{equation}
Moreover,
\begin{equation}\label{modp4}
{np \choose rp} \Big{/} {n \choose r} \equiv  1\pmod{p^4}
\end{equation}
for all $n,r$ if and only if $p$ divides the numerator of
$B_{p-3}$.
\end{thm}
\begin{rem} When $p=5$ Theorem~\ref{thm:1} does not
hold. Indeed, it's easy to see that $H_1(4)=25/12$ so $w_5=23$.
Now take $n=4$ and $r=1$. Then
$${4\cdot 5 \choose 5} \Big{/} {4\choose 1} \equiv 751
\not\equiv  1+23\cdot 4\cdot 1\cdot 3\cdot 5^3 \equiv
126\pmod{5^5}.$$
\end{rem}\begin{proof}
Clearly we may assume $n>r$. To save space we write
$H_k=H(\{1\}^k;p-1)$ throughout this proof. By equation \eqref{fx}
we have
 $${np \choose rp}= \frac{\prod_{j=n-r+1}^n f_p(jp)}
 { \prod_{l=1}^r f_p(lp) }.$$
By relation \eqref{Stirling} and Lemma \ref{genWols} we have
\begin{equation}\label{midd}
{np \choose rp}\Big{/}  {n\choose r}\equiv
 \frac{\prod_{j=n-r+1}^n\big(1-jpH_1 +j^2p^2 H_2\big)}
 {\prod_{l=1}^r\big(1-lpH_1+l^2p^2 H_2\big)} \pmod{p^5 }.
\end{equation}
Now it follows quickly from \eqref{equ:wols} and the shuffle
relation \eqref{shuf} that
\begin{equation}\label{H12}
2H_2+H(2;p-1)=H_1^2\equiv 0\pmod{p^4}.
\end{equation}
By substitution $k\to p-k$ we further can see that
\begin{alignat}{3}
2H_1 =& \, \sum_{k=1}^{p-1}\frac{p}{k(p-k)} \equiv -
\sum_{k=1}^{p-1} \frac{p}{k^2}\left(1+\frac{p}{k}+\frac{p^2}{k^2}
 \right) &\pmod{p^4}\phantom{.} \notag \\
\equiv&\,  -\big(pH(2;p-1)+p^2H(3;p-1)+p^3H(4;p-1)\big) )
&\pmod{p^4}\phantom{.} \notag \\
\equiv&\, -pH(2;p-1)  &\pmod{p^4}\phantom{.} \label{H12p} \\
\equiv&\, 2pH_2 &\pmod{p^4}. \notag
\end{alignat}
Congruence \eqref{H12p} is obtained by Lemma~\ref{genWols} (so we
indeed need the condition $p\ge 7$) while the last step follows
from \eqref{H12}. Therefore by \eqref{equ:wols} congruence
\eqref{midd} is reduced to
\begin{alignat*}{3}
{np \choose rp}\Big{/}{n\choose r} \equiv&\,
  1+\sum_{j=n-r+1}^n (j^2-j)pH_1
 -\sum_{l=1}^r  (l^2-l)pH_1  &\pmod{p^5 }\phantom{.} \\
 \equiv& \, 1+w_pnr(n-r)p^3  &\pmod{p^5 }.
\end{alignat*}
This proves congruence \eqref{main}. The last statement of the
theorem follows from \eqref{Gl} immmediately.
\end{proof}
By induction on the exponent the following corollary is obvious.
\begin{cor}\label{exp}
Let $p\ge 7$ be a prime and let $r$ and $n$ be two non-negative
integers. Then for any exponent $e\ge 1$ we have
$${np^e \choose rp^e} \Big{/} {n \choose r}\equiv 1+w_pnr(n-r)p^3 \pmod{p^5}.$$
\end{cor}

Next we consider a refined version of Theorem~\ref{b2} of Bailey.
\begin{thm} \label{thm:2}
Let $N$, $R$, $n$ and $r$ be non-negative integers and $p\ge 7$ be
a prime. If $r\le n<p$ then
\begin{equation}\label{1cong}
 {Np^3+n \choose Rp^3+r}\Big{/} \Big[{N\choose R}{n\choose r}\Big]
  \equiv 1+c(N,R,n,r;p)p^3   \pmod{p^5},
\end{equation}
where $c(N,R,n,r;p)=H_1(n)N-H_1(r)R+(w_pNR-H_1(n-r))(N-R)$. If
$n<r<p$ then
\begin{equation}\label{2cong}
{Np^3+n \choose Rp^3+r} \Big{/} {N\choose R}\equiv (-1)^{r-n+1}
\frac{N-R}{n} {r-1\choose n}^{-1} p^3\pmod{p^5}.
\end{equation}
\end{thm}

\begin{proof} Clearly we can assume that $N\ge R$.
Observe that we have a variant of $f_n(x)$ defined by \eqref{fx}:
\begin{align*}
F_n(x)=&(-1)^{n+1}f_{n+1}(-x)/x=(x+1)(x+2)\cdots(x+n)\\
=&n!\big(1+H(1;n)x+H(1,1;n)x^2+\cdots\big).
\end{align*}
First let $r\le n<p$. Then by straight-forward expansion we have
\begin{alignat*}{3}
{Np^3+n \choose Rp^3+r}\Big/
 \Big[{N \choose R} {n \choose r} \Big]  =
 & \Big[ {N p^3\choose Rp^3}\Big/ {N\choose R} \Big]
 \frac{r!(n-r)!\cdot \prod_{i=1}^n(Np^3+i)/n!}{\prod_{i=1}^r
 (Rp^3+i)\prod_{i=1}^{n-r}((N-R)p^3+i)} &\ \\
 = &  \Big[ {N p^3\choose Rp^3}\Big/ {N\choose R} \Big]
  \frac{F_n(Np^3)}{F_r(Rp^3) F_{n-r}((N-R)p^3)} &\ \\
 \equiv &\frac{(1+H_1(n)Np^3)(1+w_pNR(N-R)p^3)}{(1+H_1(r)Rp^3)(1+H_1(n-r)(N-R)p^3)}
 &\pmod{p^5}
\end{alignat*}
by Corollary \ref{exp} and the fact that $H_1(m)$ is $p$-integral
if $m<p$. Congruence \eqref{1cong} follows immediately.
\end{proof}

\begin{eg} Take $p=7$. Then the following congruence is exact
(and the term $c(\cdots)p^3$ is not needed):
 $${4\cdot 7^3+ 5 \choose 2 \cdot 7^3+2} \equiv
 {4\choose 2}{5\choose 2} \pmod{7^5}.$$
Using GP Pari and taking $1\le N,R,n,r\le 6$ we find the complete
list of nontrivial $(N,R,n,r)$ (i.e., $N\ne R$ or $n\ne r$) for
which this type of congruence holds when $p=7$: (4,2,5,2),
(4,2,5,3), (5,2,6,1), (4,2,6,3), (5,1,6,3), (5,4,6,3), (5,3,6,5).
We believe there are always such nontrivial congruences for every
prime $p\ge 7$.
\end{eg}

\begin{rem} (1) If $p=5$ then the congruence \eqref{1cong} of
Theorem~\ref{thm:2} is not true anymore. For
example, take $N=3,n=4,R=r=1$. Then $c(3,1,4,1;5)=1675/12$. So
$${3\cdot 5^3+4\choose 5^3+1}\Big{/} \Big[{3\choose 1}{4\choose 1}\Big]\equiv
2501 \not\equiv  1+c(3,1,4,1;5)5^3  \equiv 1\pmod{5^5}.$$

(2) If $p=5$ then congruence \eqref{2cong} of Theorem~\ref{thm:2}
still holds for all possible $N,R<5^5$ and $n<r<5$. I believe this
is true for all other $N$ and $R$.
\end{rem}

\section{An interesting sum related to $\zeta(1,2)$}
The last result of this note is related to the above theme and has
some independent interest. We discovered this when trying to prove
Theorem~\ref{thm:1} in the special case $r=1$ following Gardiner's
suggestion in \cite{Gardiner}. We failed but obtained this
unexpected byproduct.
\begin{prop} \label{prop:ijk}
Suppose $p$ is an odd prime. Then
$$2H(2,1;p-1)\equiv -2H(1,2;p-1)\equiv
\sum_{\substack {i+j+k=p \\ i,j,k>0}} \frac{1}{ijk}  \pmod {p}.$$
\end{prop}
\begin{proof} By the shuffle relation \eqref{shuf} and
Lemma~\ref{genWols} we have
 $$H(2,1;p-1)+H(1,2;p-1)=H(1;p-1)H(2;p-1)-H(3;p-1)\equiv 0 \pmod{p}.$$
So the first congruence is obvious. Let's prove the second.

The cases $p=3$ and $5$ can be checked easily:
\begin{alignat*}{3}
2H(1,2;2)+\sum_{\substack {i+j+k=3 \\ i,j,k>0}} \frac{1}{ijk}
=&\frac{1}{2}+1=\frac{3}{2}\equiv 0 &\pmod 3,\\
2H(1,2;4) + \sum_{\substack {i+j+k=5 \\ i,j,k>0}} \frac{1}{ijk}
=&\frac{17}{16}+\frac{7}{4}=\frac{45}{16} \equiv 0 &\pmod 5.
\end{alignat*}

Suppose now $p\ge 7$. Let's go through Gardiner's proof of
\cite[Theorem~1]{Gardiner}. Let $n>3$ be a positive integer (we
will take $n=p-1$ later). Combinatorial consideration leads us to
\begin{align}{3}
{np\choose p}
=&\sum_{i_1+\dots+i_n=p} {p\choose i_1} \cdots {p\choose i_n}\notag\\
\equiv & n+ {n\choose 2} \sum_{\substack {i+j=p\\ i,j>0}}
{p\choose i}{p\choose j}+
 {n\choose 3}\sum_{\substack {i+j+k=p \\ i,j,k>0}}
 {p\choose i}{p\choose j}{p\choose k}  & \pmod{p^4}\label{2sum}\\
\equiv &n+{n\choose 2}X+ {n\choose 3} Y  & \pmod{p^4}\notag
\end{align}
where $X$ and $Y$ are given by the two sums in \eqref{2sum}
respectively. Recall from \eqref{fx} and \eqref{Stirling}
$$f_i(x)=x(x-1)\cdots(x-i+1)=(i-1)!\sum_{j=1}^n (-1)^{i-j}H(\{1\}^{j-1};i-1)  x^j.$$
Hence
\begin{alignat*}{3}
X=\sum_{i=1}^{p-1}\left(\frac{f_i(p)}{i!}\right)^2
=&\sum_{i=1}^{p-1}\frac{p^2}{i^2} \left(
\sum_{j=1}^{i} (-1)^{i-j}H(\{1\}^{j-1};i-1) p^{j-1} \right)^2 &\  \\
\equiv & \sum_{i=1}^{p-1} \frac{p^2}{i^2}\big(1-2H(1;i-1)p \big) & \pmod {p^4}\phantom{.}\\
\equiv &  p^2 H(2;p-1) -2p^3 H(1,2;p-1)  & \pmod {p^4}.
\end{alignat*}
As for $Y$ we have $Y=0$ if $n=2$. If $n\ge 3$ then
\begin{align*}
Y=& \sum_{\substack {i+j+k=p \\ i,j,k>0}}
\frac{1}{i!j!k!}\prod_{\ga=i,j,k} \left(
    \sum_{l=1}^{\ga}(-1)^{l-\ga}S(\ga,l) p^l \right)
\equiv \sum_{\substack {i+j+k=p \\ i,j,k>0}} \frac{p^3}{ijk} \pmod
{p^4}.
\end{align*}
Putting every thing together with $n=p-1$, comparing to
\eqref{main} with $r=1$, using the fact $2H(1;p-1)\equiv
-pH(2;p-1) \pmod{p^4}$ from \eqref{H12p}, and canceling the factor
$n(n-1)/2$, we arrive at
$$
- n p^2H(2;p-1) \equiv p^2H(2;p-1) -2p^3 H(1,2;p-1) +\frac{n-2}{3}
\sum_{\substack {i+j+k=p \\ i,j,k>0}} \frac{p^3}{ijk} \pmod {p^4}.
$$
With $n=p-1$ this simplifies to
$$2 H(1,2;p-1)-\frac{p-3}{3}\sum_{\substack{i+j+k=p \\ i,j,k>0}}
 \frac{1}{ijk}\equiv H(2;p-1) \equiv 0 \pmod {p},$$
whence the second congruence in the proposition.
\end{proof}

Combining Proposition~\ref{prop:ijk} with \cite[Theorem~3.1]{wols}
we find the following corollary.
\begin{cor} \label{cor:ijk}
For any prime $p\ge 5$
$$\sum_{\substack {i+j+k=p \\ i,j,k>0}} \frac{1}{ijk}  \equiv
-2B_{p-3} \pmod {p}.$$
\end{cor}

\begin{rem} (1) We know that among all the primes $p$ less than 12 million
$p$ divides the numerator of $B_{p-3}$ only for $p=16843$ and
$p=2124679$ (see \cite{BCEMS}). However, we believe there exist
infinitely many such primes.

(2) In a recent paper Ji provides a proof of Corollary
\ref{cor:ijk} without using partial sums of MZV series (see
\cite{Ji}). In more recent preprints he \cite{J2} and
independently, Zhou and Cai \cite{zc}, generalize this to sums of
arbitrary lengths: let $p\ge 5$ be a prime and $n\leq p-2$ a
positive integer then
    $$\sum_{\substack {l_{1}+l_{2}+\cdots+l_{n}=p\\l_{1},\ldots,\, l_{n}>0}}
    \frac{1}{l_{1}l_{2}\cdots l_{n}}\equiv
   \left \{\aligned
 -(n-1)!\ B_{p-n} \quad \quad  \  \pmod{p} \quad & \text{if } 2\nmid n,\,\\
 -\frac{ n!np}{2(n+1)} B_{p-n-1}\ \ \pmod{p^2}\quad  &\text{if } 2\, |\,n.
 \endaligned\right.
$$
\end{rem}

\bigskip

\noindent {\em Email:} zhaoj@eckerd.edu


\begin{thebibliography}{9}
\bibitem{Bailey}
D.F.~Bailey, \textit{Two $p^3$ variations of Lucas' theorem,} J.
Number Theory  \textbf{35}(2) (1990), pp.~208--215. \textbf{MR:}
91f:11008.

\bibitem{Buhcom}
J.P.~Buhler, private email correspondence, Dec. 12, 2002.

\bibitem{BCEMS}
J.P.~Buhler, R.E.~Crandall, R.~Ernvall, and T.~Mets\"ankyl\"a, and
M.A.~Shokrollahi, \textit{Primes and cyclotomic invariants to 12
million,} Computational algebra and number theory (Milwaukee, WI,
1996). J. Symbolic Comput.  \textbf{31} (2001), pp.~89--96.
\textbf{MR:} 2001m:11220.

\bibitem{Gardiner}
A.~Gardiner, \textit{Four problems on prime power divisibility},
Amer. Math. Monthly \textbf{95} (1988), pp.~926--931.

\bibitem{Gl}
J.W.L.~Glaisher,  {}{On the residues of the sums of the inverse
powers of numbers in arithmetical progression}, Quarterly J. Math.
\textbf{32} (1900), 271--288.

\bibitem{HW}
G.H.~Hardy and E.M.~Wright, \textit{An introduction to the theory
of numbers}, Clarendon press, Oxford, 1980.

\bibitem{Ji}
Chun-Gang Ji, \textit{A simple proof of a curious congruence by
Zhao}, Proc. Amer. Math. Soc. \textbf{133}(2005), pp.~3469-3472.

\bibitem{J2} Chun-Gang Ji, \textit{Generalization of Wolstenholme's
Theorem}, To appear in Disc. Math.

\bibitem{K1}
G.~S.~Kazandzidis, \textit{Congruences on the binomial
coefficients,} Bull. Soc. Math. Grèce (N.S.) \textbf{9} (1968),
fasc. 1, pp.~1--12. \textbf{MR:} 42\#182.

\bibitem{K2}
G.~S.~Kazandzidis, \textit{On congruences in number-theory,} Bull.
Soc. Math. Grèce (N.S.) \textbf{10} (1969), fasc. 1, pp.~35--40.
\textbf{MR:} 43\#4753.

\bibitem{ELeh}
E.~Lehmer, \textit{On Congruences Involving Bernoulli Numbers and
the Quotients of Fermat and Wilson,} Ann. Math., 2nd Ser.,
\textbf{39} (1938), pp.~350--360.

\bibitem{R}
A.M.~Robert, \textit{A Course in p-adic Analysis}, Graduate Texts
in Mathematics, Vol.~198,  Springer, 2000.

\bibitem{RZ}
A.~Robert and M.~Zuber, \textit{The Kazandzidis supercongruences.
A simple proof and an application,} Rend. Sem. Mat. Univ. Padova
\textbf{94} (1995), 235--243. \textbf{MR:} 96m:11014.

\bibitem{wols}
J.~Zhao,  \textit{Partial sums of multiple zeta value series I:
generalizations of Wolstenholme's Theorem},
xxx.lanl.gov/abs/math.NT/0301252, v1.

\bibitem{zc}
X.~Zhou and T.~Cai, \textit{A generalization of a curious
congruence by Zhao,} To appear in  Proc. Amer. Math. Soc.

\end{thebibliography}
\end{document}